\documentclass[a4paper,12pt]{article}

\usepackage{latexsym}
\usepackage{amssymb}
\usepackage{amsthm}

\title{On the Comultiplication in Quantum Affine Algebras}
\author{Jesper Thor\'en}
\date{}

\newtheorem{definition}{Definition}[section]
\newtheorem{theorem}[definition]{Theorem}
\newtheorem{proposition}[definition]{Proposition}
\newtheorem{lemma}[definition]{Lemma}
\newtheorem{corollary}[definition]{Corollary}

\begin{document}
\maketitle

\begin{abstract}
  We express the comultiplication of the generators in Drinfelds
  second realization of the quantum affine algebra
  $U_q(\hat{\mathfrak{sl}}_{2})$, induced by the
  comultiplication of the generators in the Drinfeld-Jimbo realization of
  $U_q(\hat{\mathfrak{sl}}_{2})$ in terms of generating
  functions. Then we find explicit expressions for the
  comultiplication of the generators. 
\end{abstract}

\section*{Introduction}

In \cite{Dr2}, V.G. Drinfeld gave a new set of generators and relations for
the quantum affine algebra $U_q(\hat{\mathfrak {g}})$ (and also for the Yangian). He also gave an
isomorphism between the two realizations, but there was no proof in
that article. In \cite{Be}, J. Beck found these new generators inside
$U_q(\hat{\mathfrak {g}})$, and proved that they satisfy the
relations given by Drinfeld. He also proved that the
two realizations were isomorphic as  Hopf algebras, and an
explicit isomorphism was given in that article. Beck also gave new
formulas satisfied by the comultiplication, but no explicit
expressions for the comultiplication of the generators were
found. Another comultiplication was found by Drinfeld in an
unpublished paper, see \cite{DiFr}, \cite{Di}, but this new comultiplication has
values in an extension of $U_q(\hat{\mathfrak {g}})\otimes
U_q(\hat{\mathfrak {g}})$, and do not coincide with the
comultiplication induced by the isomorphism given by Beck in \cite{Be}. Its advantage over the
induced comultiplication is that it is much easier to work with (cf.
\cite{Ji}, \cite{Di}). 

In this paper we find the comultiplication of
the generators in Drinfelds new realization of the quantum affine algebra
$U_q(\hat{\mathfrak{sl}}_{2})$, using Beck's isomorphism (see
\cite{Th} for the corresponding case for the Yangian
associated to $\hat{\mathfrak{sl}}_{2}$). 

The paper is organized as follows: Section 1 contains the definition
of the algebras $U_q(\hat{\mathfrak{sl}}_{2})$ and $\hat{U}_q.$ We
define certain generating functions, $X_0^{\pm}(z), ~Y_0^{\pm}(z),$ and some
further relations involving the generators of $\hat{U}_q$ are found.

In Section 2, the main theorem is Theorem \ref{theorem5}, where we give
relations satisfied by the comultiplication of the generators of the
algebra $\hat{U}_q$. The result is
written in terms of the generating functions defined in Section 1.

Finally, in Section 3, we first study powers of the generating functions
$X_0^+(z),$ and $Y_0^+(z).$ One important result of this section is Theorem
\ref{theorem6}, which shows us how to express such a power as a sum of
PBW-ordered monomials. Then we deduce, from this theorem and from
Theorem \ref{theorem5}, explicit formulae for the comultiplication of
the generators.

{\em Acknowledgement. } I thank Prof. A. Meurman for helpful
discussions.

\section{Relations in $\hat{U}_q$.}

In this section we define the isomorphic  algebras $U_q(\hat{\mathfrak{sl}}_{2})$
and $\hat{U}_q.$ Then we find and prove some new relations in
$\hat{U}_q,$ and express them in terms of generating functions.

Let $q$ be an indeterminate. For integers $n,r, ~n \geq r \geq 0,$
define

$$[n]=\frac{q^n-q^{-n}}{q-q^{-1}} \in \mathbb{Z} [q,q^{-1}],
  \quad [n]! = [n][n-1]\cdots[1],$$
and
$$\left[ {n \atop r} \right] = \frac{[n]!}{[n-r]![r]!}\in \mathbb{Z} [q,q^{-1}].$$
It is also convenient to define $[0]!=1.$

The following definition is due to Drinfeld \cite{Dr1} and Jimbo
\cite{J} (note that, as in \cite{Be}, we have added the square root of
$k_0k_1$).

\begin{definition}
  The Hopf algebra $U_q(\hat{\mathfrak{sl}}_{2})$ is an associative
  algebra over $\mathbb{C}(q)$ with generators $e_i^\pm , k_i^{\pm 1},
  (i \in \{0,1\}), c^{\pm 1/2},$ and with defining relations
  \begin{eqnarray*}
    (c^{\pm 1/2})^2 = (k_0k_1)^{\pm 1}, && c^{1/2}~ {\textrm {is~ central}}, \\
    k_ik_i^{-1}=k_i^{-1}k_i=1,&& k_ik_j=k_jk_i, \\
    k_ie_i^\pm k_i^{-1}=q^{\pm 2}e_i^\pm, &&
    k_ie_j^\pm k_i^{-1}=q^{\mp 2}e_j^\pm , ~i\neq j, \\
    {} [e_i^+,e_j^-] &=& \delta_{ij} \frac{k_i-k_i^{-1}}{q-q^{-1}}, \\
    \sum_{r=0}^3 (-1)^r \left[{3 \atop r} \right]
    (e_i^\pm)^re_j^\pm(e_i^\pm)^{3-r} &=& 0, ~i\neq j. 
  \end{eqnarray*}
  
  The comultiplication of $U_q(\hat{\mathfrak{sl}}_{2})$ is given by
  \begin{eqnarray*}
    \Delta (e_i^+) &=& e_i^+ \otimes k_i + 1 \otimes e_i^+,\\
    \Delta (e_i^-) &=& e_i^- \otimes 1 + k_i^{-1} \otimes e_i^-,\\
    \Delta (k_i) &=& k_i \otimes k_i, \\
    \Delta (c^{\pm 1/2}) &=& c^{\pm 1/2}\otimes c^{\pm 1/2}.
  \end{eqnarray*}
\end{definition}

In \cite{Dr2}, V.G. Drinfeld gave a new set of generators and relations for
$U_q(\hat{\mathfrak{sl}}_{2})$. We formulate this new version of
$U_q(\hat{\mathfrak{sl}}_{2})$ in the following proposition. See J. Beck,
\cite{Be}, for a proof.

Let $z$ be an indeterminate.

\begin{proposition} \label{proposition1}
  Let $\hat{U}_q$ be the $\mathbb{C}(q)$-algebra generated by
  $$x_n,~ y_n, ~ (n\in \mathbb{Z}),~ h_k, ~ (k \in
  \mathbb{Z}\smallsetminus \{ 0 \}),~ K^{\pm
    1},~ c^{\pm 1/2},$$
  and with the following defining relations
  \begin{eqnarray}
    KK^{-1}=1, && K^{-1}K=1, \label{prop1.1} \\
    c^{1/2}c^{-1/2}=c^{-1/2}c^{1/2}=1, && c^{1/2}~ {\textrm {is~ central}}, \label{prop1.2} \\
    {} [h_m,h_n] &=& \delta_{m,-n} \frac{[2m]}{m}
    \frac{c^m-c^{-m}}{q-q^{-1}}, \label{prop1.3} \\
    Kh_m &=& h_mK, \label{prop1.4} \\
    Kx_mK^{-1}=q^2x_m, & &  Ky_mK^{-1}=q^{-2}y_m, \label{prop1.5} \\
    {} [h_m,x_n] &=& \frac{[2m]}{m} c^{-|m|/2} x_{m+n},
    \label{prop1.6} \\
    {} [h_m,y_n] &=& -\frac{[2m]}{m} c^{|m|/2} y_{m+n},
    \label{prop1.7} \\
    x_{m+1}x_n-q^2x_nx_{m+1} &=& q^2x_mx_{n+1}-x_{n+1}x_m,
    \label{prop1.8} \\
    y_{m+1}y_n-q^{-2}y_ny_{m+1} &=& q^{-2}y_my_{n+1}-y_{n+1}y_m,
    \label{prop1.9} \\
    {} [x_m,y_n] &=&
    \frac{c^{(m-n)/2}\psi_{m+n}-c^{-(m-n)/2}\phi_{m+n}}{q-q^{-1}},
    \label{prop1.10} 
  \end{eqnarray}  
  where we have set $\psi_{-m} = 0,~ \phi_m = 0,~ m >0$, and where $\psi_m, $  and
  $\phi_{-m},~ m \in \mathbb{N}$ are given by
  $$\sum_{m=0}^\infty \psi_m z^{-m} = K \exp \Big(
  (q-q^{-1})\sum_{m=1}^\infty h_m z^{-m} \Big),$$
  and
  $$\sum_{m=0}^\infty \phi_{-m} z^{m} = K^{-1} \exp \Big(
  -(q-q^{-1})\sum_{m=1}^\infty h_{-m} z^{m} \Big),$$

    Then $U_q(\hat{\mathfrak{sl}}_{2})$ is isomorphic to $\hat{U}_q$,
    and the isomorphism is given by
    
    $$c^{\pm 1/2}\mapsto c^{\pm 1/2}, \quad k_0 \mapsto cK^{-1}, \quad k_1 \mapsto K,$$
    $$e_1^+ \mapsto x_0, \quad  e_1^- \mapsto y_0,$$
    $$e_0^+ \mapsto y_1K^{-1}, \quad  e_0^- \mapsto Kx_{-1}.$$
    
  \end{proposition}

Let $\mathcal{X}$ be the $\mathbb{C}(q)$-subalgebra of $\hat{U}_q$ generated by the
$$x_n, ~(n \in \mathbb{Z}),~K^{\pm 1},~c^{\pm 1/2},~h_k, ~(k\in
\mathbb{Z}\smallsetminus\{ 0\}), $$
and let $\mathcal{Y}$ be the $\mathbb{C}(q)$-subalgebra of $\hat{U}_q$
generated by the
$$y_n, ~(n \in \mathbb{Z}),~K^{\pm 1},~c^{\pm 1/2},~h_k,
{}~(k\in \mathbb{Z}\smallsetminus\{ 0\}). $$
Then it is easy to prove the
following proposition, using Proposition \ref{proposition1}.

\begin{proposition} \label{proposition2}
  There exist $\mathbb{C}(q)$-algebra automorphisms
  $S:\mathcal{X}\to \mathcal{X}$ and $T:\mathcal{Y}\to
  \mathcal{Y},$ given by
  $$S:x_n\mapsto x_{n+1}, ~h_k\mapsto h_k, ~c^{\pm 1/2}\mapsto c^{\pm
    1/2}, ~K\mapsto K,$$
  and
  $$T:y_n\mapsto y_{n+1}, ~h_k\mapsto h_k, ~c^{\pm 1/2}\mapsto c^{\pm
    1/2}, ~K\mapsto K,$$
  for all $n\in \mathbb{Z}, ~ k\in \mathbb{Z}\smallsetminus \{ 0\}.$
  
\end{proposition}  
   
As in \cite{Ch}, \cite{ChPr} we have the following

\begin{proposition} \label{proposition3}
  (a). There exists a $\mathbb{C}$-algebra antiautomorphism
  $\alpha:\hat{U}_q \to \hat{U}_q$ such that
  for all $n\in \mathbb{Z}, ~ k\in \mathbb{Z}\smallsetminus \{ 0\},$
  $$\alpha:q\mapsto q^{-1}, ~c^{1/2}\mapsto c^{-1/2}, ~x_n\mapsto y_{-n},
  {}~y_n\mapsto x_{-n},$$
  $$\psi_n\mapsto \phi_{-n}, ~\phi_n\mapsto \psi_{-n}, ~h_k\mapsto h_{-k}.$$
  Moreover,
  $\alpha\otimes\alpha\circ\tau\circ\Delta=\Delta\circ\alpha$, where
  $\tau$ is the flip map given by $a\otimes b \mapsto b\otimes a$ on
  $\hat{U}_q\otimes \hat{U}_q.$ \\
  (b). There exists a $\mathbb{C}(q)$-algebra antiautomorphism
  $\beta:\hat{U}_q \to \hat{U}_q$ such that
  for all $n\in \mathbb{Z}, ~ k\in \mathbb{Z}\smallsetminus \{ 0\},$
  $$\beta:c^{1/2}\mapsto c^{-1/2}, ~x_n\mapsto y_n,
  {}~y_n\mapsto x_n,$$
  $$\psi_n\mapsto \psi_n, ~\phi_n\mapsto \phi_n, ~h_k\mapsto h_k.$$
\end{proposition}

\begin{proof}
  This is easy to check, using Proposition \ref{proposition1}.
\end{proof}
  
For $z$ an indeterminate, define the generating functions
\begin{eqnarray}
  X^+(z)=\sum_{k \geq 1} x_kz^{-k}, && X_0^+(z)=x_0+X^+(z),
  \label{defX+} \\
  X^-(z)=\sum_{k \geq 1} x_{-k}z^{k}, && X_0^-(z)=x_0+X^-(z),
  \label{defX-} \\
  Y^+(z)=\sum_{k \geq 1} y_kz^{-k}, && Y_0^+(z)=y_0+Y^+(z),
  \label{defY+} \\
  Y^-(z)=\sum_{k \geq 1} y_{-k}z^{k}, && Y_0^-(z)=y_0+Y^-(z).
  \label{defY-} 
\end{eqnarray}
We set
\begin{equation}
  X(z)=X_0^+(z)+X^-(z), \quad Y(z)=Y_0^+(z)+Y^-(z). \label{defXY}
\end{equation}
We also define
\begin{eqnarray}
  \Psi (z)=\sum_{k \geq 0}\psi_kz^{-k}, && \Phi (z)=\sum_{k \geq
    0}\phi_{-k}z^k. \label{defFiPsi}
\end{eqnarray}

We will extend the maps $\alpha$ and $\beta$ linearly to $\hat{U}_q
[[z,z^{-1}]]$, so that we can consider relations of the type
$$\alpha (Y^+(q^2cz^{-1}))=X^-(q^2cz),$$
$$\alpha (\Psi (c^{3/2}z^{-1}))=\Phi (c^{3/2}z).$$

In Section 2 we will express the comultiplication of the generators $x_k,
{}~y_k, {}~\psi_k, {}~\phi_k,$ in terms of the series above, and for that we
will need some further relations in the algebra $\hat{U}_q.$ We state
them in the following lemma.

\begin{lemma} \label{lemma4}
  For $n \in \mathbb{Z}_+,$ we have
  \begin{eqnarray}
    {}[h_1,X_0^+(z)^n] &=& [2](\sum_{k=0}^{n-1}q^{-2k})
    zc^{-1/2}(X_0^+(z)^n-x_0X_0^+(z)^{n-1}); \label{lemma4.1} \\
    {}[h_1,Y^+(z)^n] &=& -[2](\sum_{k=0}^{n-1}q^{-2k})
    c^{1/2}(zY^+(z)^n-Y^+(z)^{n-1}y_1); \label{lemma4.2} \\
    X_0^+(z)^nx_0 &=& q^{-2n}x_0X_0^+(z)^n+
    (1-q^{-2n})X_0^+(z)^{n+1}; \label{lemma4.3} \\
    X^+(z)^nx_0 &=& q^{2n}x_0X^+(z)^n-
    (1-q^{2n})X^+(z)^{n+1}; \label{lemma4.4} \\
    y_1Y^+(z)^n &=& q^{-2n}Y^+(z)^ny_1+
    (1-q^{-2n})zY^+(z)^{n+1}; \label{lemma4.5} \\
    y_0Y^+(z)^n &=& q^{2n}Y^+(z)^ny_0-
    (1-q^{2n})Y^+(z)^{n+1}; \label{lemma4.6} \\
    x_0\Psi (z) &=& q^2\Psi (z)x_0-
    (q^2-q^{-2})\Psi (z)X_0^+(q^2c^{1/2}z); \label{lemma4.7} \\
    \Psi (z)y_1 &=& q^2y_1\Psi (z)+
    (1-q^{4})zc^{-1/2}Y^+(q^2c^{-1/2}z)\Psi (z); \label{lemma4.8} \\
    {}[X_0^+(z)^n,y_0] &=& (q-q^{-1})^{-1}q^{n-1}[n]
    \big( q^{-2n+2}\Psi (c^{-1/2}z)X_0^+(q^2z)^{n-1} \nonumber \\
    && -K^{-1}X_0^+(z)^{n-1}\big); \label{lemma4.9} \\
    {}[Y^+(z),x_0] &=& -(q-q^{-1})^{-1} (\Psi (c^{1/2}z)-K);
    \label{lemma4.10} \\
    {}[X_0^+(z),y_1] &=& (q-q^{-1})^{-1} zc^{-1}(\Psi (c^{-1/2}z)-K).
    \label{lemma4.11}
  \end{eqnarray}
\end{lemma}

\begin{proof}
  We will prove (\ref{lemma4.3}) first.
  Write
  $$[x,y]_{q^2}=xy-q^2yx.$$
  Suppose that $n=1.$ Proposition \ref{proposition1} (\ref{prop1.8})
  gives the relation
  $$[x_{k+1},x_l]_{q^2}=q^2[x_k,x_{l+1}]_{q^2}+(q^4-1)x_{l+1}x_k.$$
  Multiplying both sides of this relation with $z^{-k-l-1}$ and
  summing over $k, l \in \mathbb{N}$ gives
  $$[X_0^+(z)-x_0,X_0^+(z)]_{q^2}=q^2[X_0^+(z),X_0^+(z)-x_0]_{q^2}
  +(q^4-1)(X_0^+(z)^2-x_0X_0^+(z)),$$
  which is equivalent to
  $$[x_0,X_0^+(z)]_{q^2}=(1-q^2)X_0^+(z)^2,$$
  so that
  $$X_0^+(z)x_0=q^{-2}x_0X_0^+(z)+(1-q^{-2})X_0^+(z)^2.$$
  So (\ref{lemma4.3}) is true for $n=1$. Now suppose, by induction,
  that (\ref{lemma4.3}) is true for $n-1\geq 1.$ Then we have
  \begin{eqnarray*}
    X_0^+(z)^nx_0 &=& X_0^+(z)(X_0^+(z)^{n-1}x_0) \\
    &=& X_0^+(z)\Big( q^{-2n+2}x_0X_0^+(z)^{n-1}+
    (1-q^{-2n+2})X_0^+(z)^{n}\Big) \\
    &=& q^{-2n}x_0X_0^+(z)^{n}+(q^{-2n+2}-q^{-2n})X_0^+(z)^{n+1}
    +(1-q^{-2n+2})X_0^+(z)^{n+1} \\
    &=& q^{-2n}x_0X_0^+(z)^{n}+(1-q^{-2n})X_0^+(z)^{n+1}.
  \end{eqnarray*}
  So (\ref{lemma4.3}) holds, by induction, for $n\geq 1.$ Now,
  (\ref{lemma4.4}) can be proved similarly, and (\ref{lemma4.5})
  follows after applying the antiautomorphism $\beta$ from Proposition
  \ref{proposition3} on (\ref{lemma4.3}), and then the automorphism
  $T$ from Proposition \ref{proposition2} on the image of $\beta.$ For
  (\ref{lemma4.6}), use $\beta$ on (\ref{lemma4.4}).

  We can now prove (\ref{lemma4.1}) by induction over $n.$ First
  consider the case $n=1.$ By Proposition \ref{proposition1}
  (\ref{prop1.6}) we have that
  $$[h_1,\sum_{k\geq 0} x_kz^{-k}]=[2]c^{-1/2}\sum_{k\geq
    0}x_{k+1}z^{-k}.$$
  So
  $$[h_1,X_0^+(z)]=[2]c^{-1/2}z(X_0^+(z)-x_0),$$
  which proves (\ref{lemma4.1}) when $n=1.$ Now suppose
  (\ref{lemma4.1}) is true for $n-1\geq 1.$ Then, by induction,
  \begin{eqnarray*}
    [h_1,X_0^+(z)^n] &=& [h_1,X_0^+(z)]X_0^+(z)^{n-1}
    +X_0^+(z)[h_1,X_0^+(z)^{n-1}] \\
    &=& [2]c^{-1/2}z\Big( (X_0^+(z)^n-x_0X_0^+(z)^{n-1}) \\
    && +(\sum_{k=0}^{n-2}q^{-2k})(X_0^+(z)^n -
    X_0^+(z)x_0X_0^+(z)^{n-2})\Big).
  \end{eqnarray*}
  By (\ref{lemma4.3}) we get
  \begin{eqnarray*}
    [h_1,X_0^+(z)^n] &=& [2]zc^{-1/2}\Big(
    X_0^+(z)^n-x_0X_0^+(z)^{n-1} \\
    && +(\sum_{k=0}^{n-2}q^{-2k})(X_0^+(z)^n-q^{-2}x_0X_0^+(z)^{n-1} \\
    && -(1-q^{-2})X_0^+(z)^n)\Big) \\
    &=& [2]zc^{-1/2}\Big( (1+\sum_{k=0}^{n-2}q^{-2k-2})X_0^+(z)^n \\
    && -(1+\sum_{k=0}^{n-2}q^{-2k-2})x_0X_0^+(z)^{n-1}\Big) \\
    &=& [2](\sum_{k=0}^{n-1}q^{-2k})zc^{-1/2} \Big(
    X_0^+(z)^n-x_0X_0^+(z)^{n-1}\Big).
  \end{eqnarray*}
  The relation (\ref{lemma4.1}) follows by induction. Now
  (\ref{lemma4.2}) follows by applying $\beta$ on (\ref{lemma4.1}).

  We now intend to prove (\ref{lemma4.7}). Let $w$ be another
  indeterminate, commuting with $z$. Recall the notation (\ref{defXY}). As in \cite{Dr2}
  we have the relation
  \begin{equation}
    X(w)\Psi(z)=\frac{zc^{1/2}-q^2w}{q^2zc^{1/2}-w}\Psi(z)X(w), \label{drinrel}
  \end{equation}
  where the rational function in the right hand side of
  (\ref{drinrel}) is expanded into a series with nonnegative powers of
  $w$. Note that the relation (\ref{drinrel}) follows from the
  relations in Proposition \ref{proposition1}.

  Thus, we get that
  \begin{eqnarray*}
    X(w)\Psi(z) &=& \Psi(z) (q^{-2}-wz^{-1}c^{-1/2}) \sum_{k\in
      \mathbb{Z}}\sum_{m\geq 0} q^{-2m}c^{-m/2}x_kw^{m-k}z^{-m} \\
    &=& \Psi(z)\sum_{k\in \mathbb{Z}}\sum_{m\geq 0}
    q^{-2m-2}c^{-m/2}x_kw^{m-k}z^{-m} \\
    && -\Psi(z)\sum_{k\in \mathbb{Z}}\sum_{m\geq 0}
    q^{-2m}c^{-m/2-1/2}x_kw^{m-k+1}z^{-m-1}.
  \end{eqnarray*}
  The coefficient of $w^0$ is
  \begin{eqnarray*}
    x_0\Psi(z) &=& \Psi(z)\sum_{m\geq 0}q^{-2m-2}c^{-m/2}x_mz^{-m} 
    -\Psi(z) \sum_{m\geq 0}q^{-2m}c^{-m/2-1/2}x_{m+1}z^{-m-1} \\
    &=& q^{-2}\Psi(z)X_0^+(q^2c^{1/2}z) 
    -q^2\Psi(z)\Big( X_0^+(q^2c^{1/2}z)-x_0\Big) \\
    &=& q^2\Psi(z)x_0-(q^2-q^{-2})\Psi(z)X_0^+(q^2c^{1/2}z),
  \end{eqnarray*}
  and this proves (\ref{lemma4.7}). Now (\ref{lemma4.8}) follows from
  (\ref{lemma4.7}) by first applying $\beta,$ then $T.$

  We will now prove (\ref{lemma4.9}) by induction over $n\geq 1.$ By
  Proposition \ref{proposition1} (\ref{prop1.10}), we have
  $$[\sum_{k\geq 0}x_kz^{-k},y_0]=(q-q^{-1})^{-1}(\sum_{k\geq 0}
  c^{k/2}\psi_kz^{-k}-K^{-1}),$$
  that is,
  $$[X_0^+(z),y_0]=(q-q^{-1})^{-1}(\Psi(c^{-1/2}z)-K^{-1}),$$
  which is (\ref{lemma4.9}) when $n=1.$ Now suppose, by induction,
  that (\ref{lemma4.9}) is true for some $n-1\geq 1.$ Then, by using
  (\ref{lemma4.3}), (\ref{lemma4.7}), and the induction hypothesis, we get 
  \begin{eqnarray*}
    \lefteqn{(1-q^{-2n+2})[X_0^+(z)^n,y_0]} \\
    &=& [X_0^+(z)^{n-1}x_0,y_0]
    -q^{-2n+2}[x_0X_0^+(z)^{n-1},y_0] \\
    &=& [X_0^+(z)^{n-1},y_0]x_0+X_0^+(z)^{n-1}[x_0,y_0] \\
    &&
    -q^{-2n+2}[x_0,y_0]X_0^+(z)^{n-1}-q^{-2n+2}x_0[X_0^+(z)^{n-1},y_0]
    \\
    &=& (q-q^{-1})^{-1}\Bigg(
    q^{n-2}[n-1]\Big( q^{-2n+4}\Psi(c^{-1/2}z)X_0^+(q^2z)^{n-2}
    -K^{-1}X_0^+(z)^{n-2}\Big) x_0 \\
    && +X_0^+(z)^{n-1}(K-K^{-1})-q^{-2n+2}(K-K^{-1})X_0^+(z)^{n-1} \\
    && -q^{-2n+2}q^{n-2}[n-1]x_0\Big( q^{-2n+4}
    \Psi(c^{-1/2}z)X_0^+(q^2z)^{n-2} -K^{-1}X_0^+(z)^{n-2}\Big) \Bigg) \\
    &=& (q-q^{-1})^{-1}\Big( q^{-n+2}[n-1]\Psi(c^{-1/2}z)
    X_0^+(q^2z)^{n-2}x_0 \\
    && -q^{n-2}[n-1]K^{-1}X_0^+(z)^{n-2}x_0 
    +(q^{-2n+2}-q^{2n-2})K^{-1}X_0^+(z)^{n-1} \\
    && -q^{-3n+4}[n-1]x_0\Psi(c^{-1/2}z)X_0^+(q^2z)^{n-2} 
    +q^{-n+2}[n-1]K^{-1}x_0X_0^+(z)^{n-2} \Big) \\
    &=& (q-q^{-1})^{-1}\Bigg( q^{-n+2}[n-1]\Psi(c^{-1/2}z)\Big(
    q^{-2n+4}x_0 X_0^+(q^2z)^{n-2} \\
    &&+(1-q^{-2n+4})X_0^+(q^2z)^{n-1}\Big)
    -q^{n-2}[n-1]K^{-1}\Big( q^{-2n+4}x_0 X_0^+(z)^{n-2} \\
    && +(1-q^{-2n+4})X_0^+(z)^{n-1}\Big) 
    +(q^{-2n+2}-q^{2n-2})K^{-1}X_0^+(z)^{n-1} \\
    && -q^{-3n+4}[n-1]\Big( q^2\Psi(c^{-1/2}z)x_0 
    -(q^2-q^{-2})\Psi(c^{-1/2}z)X_0^+(q^2z)\Big) X_0^+(q^2z)^{n-2} \\
    && +q^{-n+2}[n-1]K^{-1}x_0X_0^+(z)^{n-2}\Bigg) \\
    &=& (q-q^{-1})^{-1}\Bigg( \Big( q^{-n+2}(1-q^{-2n+4})[n-1] \\
    && +q^{-3n+4}(q^2-q^{-2})[n-1]\Big) \Psi(c^{-1/2}z)X_0^+(q^2z)^{n-1}
    \\
    && +\Big( -q^{n-2}(1-q^{-2n+4})[n-1]+(q^{-2n+2}-q^{2n-2})\Big)
    K^{-1}X_0^+(z)^{n-1}\Bigg) \\
    &=& (q-q^{-1})^{-1}\Big(
    (q^{-n+2}-q^{-3n+2})[n-1]\Psi(c^{-1/2}z)X_0^+(q^2z)^{n-1} \\
    && -q^{n-1}(1-q^{-2n+2})[n] K^{-1}X_0^+(z)^{n-1}\Big) \\
    &=& (q-q^{-1})^{-1}\Big(
    q^2q^{-2n}(q^n-q^{-n})[n-1]\Psi(c^{-1/2}z)X_0^+(q^2z)^{n-1} \\
    && -q^{n-1}(1-q^{-2n+2})[n] K^{-1}X_0^+(z)^{n-1}\Big) \\
    &=& (q-q^{-1})^{-1}q^{n-1}(1-q^{-2n+2})[n]\Big(
    q^{-2n+2}\Psi(c^{-1/2}z)X_0^+(q^2z)^{n-1} 
    -K^{-1}X_0^+(z)^{n-1}\Big) ,
  \end{eqnarray*}
  and (\ref{lemma4.9}) follows.

  Finally, (\ref{lemma4.10}) and (\ref{lemma4.11}) follows from
  $$[\sum_{k\geq 1}y_kz^{-k},x_0]=-(q-q^{-1})^{-1}\sum_{k\geq
    1}c^{-k/2}\psi_kz^{-k},$$
  and
  $$[\sum_{k\geq 0}x_kz^{-k},y_1]=(q-q^{-1})^{-1}zc^{-1}\sum_{k\geq
    1}c^{k/2}\psi_kz^{-k},$$
  respectively.

  This completes the proof of the lemma.
\end{proof}

\section{The comultiplication}

The main goal of this section is to formulate and prove Theorem
\ref{theorem5}. For this we will use the generating functions defined
in (\ref{defX+})-(\ref{defFiPsi}).

Using Proposition \ref{proposition1}, we find that
\begin{eqnarray*}
  \Delta (c^{\pm 1/2}) &=& c^{\pm 1/2}\otimes c^{\pm 1/2}, \\
  \Delta(K^{\pm 1}) &=&  K^{\pm 1}\otimes K^{\pm 1} \\
  \Delta(x_0) &=& x_0\otimes K + 1\otimes x_0, \\
  \Delta(y_0) &=& y_0\otimes 1 + K^{-1}\otimes y_0, \\
  \Delta(x_{-1}) &=& x_{-1}\otimes K^{-1} + c^{-1}\otimes x_{-1}, \\
  \Delta(y_1) &=& y_1\otimes c + K\otimes y_1.
\end{eqnarray*}
Since
$$h_1 = (q-q^{-1})^{-1}K^{-1}\psi_1,$$
by Proposition \ref{proposition1}, we get that
\begin{eqnarray*}
  \Delta(h_1) &=& \Delta(c^{1/2}K^{-1})[\Delta(x_0),\Delta(y_1)] \\
  &=& h_1\otimes c^{3/2}+c^{1/2}\otimes h_1
  -(q^2-q^{-2})c^{1/2}x_0\otimes c^{1/2}y_1.
\end{eqnarray*}

We can now prove

\begin{theorem} \label{theorem5}
  With notations as in (\ref{defX+})-(\ref{defFiPsi}), the
  comultiplication in $\hat{U}_q$ satisfies
  \begin{eqnarray}
    \lefteqn{\sum_{k\geq 0} \Delta(x_k)(zc\otimes c^2)^{-k}} \nonumber
    \\
    &=& 1\otimes X_0^+(c^2z)+ \sum_{n\geq 0}(-q(q-q^{-1})^2)^n
    X_0^+(cz)^{n+1}\otimes Y^+(q^2cz)^n\Psi(c^{3/2}z); \label{thm5.1}
    \\
    \lefteqn{\sum_{k\geq 1} \Delta(x_{-k})(zc\otimes 1)^k} \nonumber
    \\
    &=& 1\otimes X^-(z)+ \sum_{n\geq 0}(-q(q-q^{-1})^2)^n
    X^-(cz)^{n+1}\otimes Y_0^-(q^2cz)^n\Phi(c^{1/2}z); \label{thm5.2}
    \\
    \lefteqn{\sum_{k\geq 1} \Delta(y_{k})(z\otimes c)^{-k}} \nonumber
    \\
    &=& Y^+(z)\otimes 1+ \sum_{n\geq 0}(-q^{-1}(q-q^{-1})^2)^n
    \Psi(c^{1/2}z)X_0^+(q^2cz)^{n}\otimes Y^+(cz)^{n+1}; \label{thm5.3}
    \\
    \lefteqn{\sum_{k\geq 0} \Delta(y_{-k})(zc^2\otimes c)^k} \nonumber
    \\
    &=& Y_0^-(c^2z)\otimes 1+ \sum_{n\geq 0}(-q^{-1}(q-q^{-1})^2)^n
    \Phi(c^{3/2}z)X^-(q^2cz)^{n}\otimes Y_0^-(cz)^{n+1}; \label{thm5.4}
    \\
    \lefteqn{\sum_{k\geq 0} \Delta(\psi_{k})(zc^{1/2}\otimes
      c^{3/2})^{-k}} \nonumber \\
    &=& \sum_{n\geq 0}(-1)^n(q-q^{-1})^{2n}[n+1]
    \Psi(c^{1/2}z)X_0^+(q^2cz)^{n}\otimes
    Y^+(q^2cz)^{n}\Psi(c^{3/2}z); \label{thm5.5} \\
    \lefteqn{\sum_{k\geq 0} \Delta(\phi_{-k})(zc^{3/2}\otimes
      c^{1/2})^{k}} \nonumber \\
    &=& \sum_{n\geq 0}(-1)^n(q-q^{-1})^{2n}[n+1]
    \Phi(c^{3/2}z)X^-(q^2cz)^{n}\otimes
    Y_0^-(q^2cz)^{n}\Phi(c^{1/2}z). \label{thm5.6}
  \end{eqnarray}
\end{theorem}  

\begin{proof}
  We start to prove (\ref{thm5.1}). From Proposition
  \ref{proposition1} we get the relation
  $$ [h_1,x_k]=[2]c^{-1/2}x_{k+1}.$$
  This relation implies the functional equation
  \begin{eqnarray*}
    \lefteqn{[\Delta(h_1),\sum_{k\geq 0} \Delta(x_k)(zc\otimes
      c^2)^{-k}]} \\
    &=& [2]z(c^{1/2}\otimes c^{3/2})\Big( \sum_{k\geq 0}
    \Delta(x_k)(zc\otimes c^2)^{-k}-\Delta(x_0)\Big) ,
  \end{eqnarray*}
  i.e., the left hand side of (\ref{thm5.1}) is a solution of the
  equation
  \begin{equation}
    [\Delta(h_1),F(z)]=[2]z(c^{1/2}\otimes c^{3/2})(F(z)-\Delta(x_0)),\label{funct1}
  \end{equation}
  where
  $$F(z)=\sum_{k\geq 0}f_kz^{-k},~f_0=\Delta(x_0).$$
  By uniqueness of solutions of functional equations with initial
  values, it is enough to show that the right hand side of
  (\ref{thm5.1}) also is a solution of (\ref{funct1}).
  
  First we calculate the two expressions
  $$[\Delta(h_1), 1\otimes X_0^+(c^2z)+X_0^+(cz)\otimes
  \Psi(c^{3/2}z)],$$
  and
  $$[\Delta(h_1), \sum_{n\geq 1}(-q(q-q^{-1})^2)^n
  X_0^+(cz)^{n+1}\otimes Y^+(q^2cz)^n\Psi(c^{3/2}z)],$$
  where
  $$\Delta(h_1) = h_1\otimes c^{3/2}+c^{1/2}\otimes h_1
  -(q^2-q^{-2})c^{1/2}x_0\otimes c^{1/2}y_1.$$
  After that, we add the results and see that (\ref{funct1}) is
  satisfied.

  Now,
  \begin{eqnarray*}
    \lefteqn{[h_1\otimes c^{3/2},1\otimes X_0^+(c^2z)+X_0^+(cz)\otimes
      \Psi(c^{3/2}z)]} \\
    &=& 0+[h_1,X_0^+(cz)]\otimes c^{3/2}\Psi(c^{3/2}z) \\
    &=& [2](zc^{1/2}\otimes c^{3/2})(X_0^+(cz)-x_0)\otimes
    \Psi(c^{3/2}z),
  \end{eqnarray*}
  by Lemma \ref{lemma4} (\ref{lemma4.1}) (with $z$ replaced by $zc$),
  and
  \begin{eqnarray*}
    \lefteqn{[c^{1/2}\otimes h_1,1\otimes X_0^+(c^2z)+X_0^+(cz)\otimes
      \Psi(c^{3/2}z)]} \\
    &=& [2](zc^{1/2}\otimes c^{3/2})1\otimes (X_0^+(c^2z)-x_0),
  \end{eqnarray*}
  again by Lemma \ref{lemma4} (\ref{lemma4.1}).

  Further, by Lemma \ref{lemma4} (\ref{lemma4.11}), (\ref{lemma4.3})
  and (\ref{lemma4.8}),
  \begin{eqnarray*}
    \lefteqn{[x_0\otimes y_1, 1\otimes X_0^+(c^2z)+X_0^+(cz)\otimes
      \Psi(c^{3/2}z)]} \\
    &=& -(q-q^{-1})^{-1}zx_0\otimes c(\Psi(c^{3/2}z)-K)
    +x_0X_0^+(cz)\otimes y_1\Psi(c^{3/2}z) \\
    && -\Big( q^{-2}x_0X_0^+(cz)+(1-q^{-2})X_0^+(cz)^2\Big) \\
    && \otimes
    \Big( q^2y_1\Psi(c^{3/2}z)+(1-q^4)zcY^+(q^2cz)\Psi(c^{3/2}z)\Big)
    \\
    &=& -(q-q^{-1})^{-1}zx_0\otimes
    c\Psi(c^{3/2}z)+(q-q^{-1})^{-1}zx_0\otimes cK \\
    && -q^{-2}(1-q^4)zx_0X_0^+(cz)\otimes cY^+(q^2cz)\Psi(c^{3/2}z) \\
    && -q^2(1-q^{-2})X_0^+(cz)^2\otimes y_1\Psi(c^{3/2}z) \\
    && -(1-q^{-2})(1-q^4)zX_0^+(cz)^2\otimes cY^+(q^2cz)\Psi(c^{3/2}z).
  \end{eqnarray*}
  It follows that
  \begin{eqnarray}
    \lefteqn{[h_1\otimes c^{3/2}+c^{1/2}\otimes
      h_1} \nonumber \\
    \lefteqn{~~-(q^2-q^{-2})c^{1/2} x_0\otimes c^{1/2}y_1, 
    1\otimes X_0^+(c^2z)+X_0^+(cz)\otimes
      \Psi(c^{3/2}z)]} \nonumber \\
    &=& [2]z(c^{1/2}\otimes c^{3/2})\Big( 1\otimes
    X_0^+(c^2z)+X_0^+(cz)\otimes \Psi(c^{3/2}z) -x_0\otimes K-1\otimes
    x_0\Big) \nonumber \\
    && -[2]z(c^{1/2}\otimes
    c^{3/2})(q-q^{-1})^2(q+q^{-1})x_0X_0^+(cz)\otimes
    Y^+(q^2cz)\Psi(c^{3/2}z) \nonumber \\
    && +[2](c^{1/2}\otimes
    c^{1/2})q(q-q^{-1})^2X_0^+(cz)^2\otimes y_1\Psi(c^{3/2}z)
    \nonumber \\
    && -[2]z(c^{1/2}\otimes
    c^{3/2})(q-q^{-1})^3(1+q^2)X_0^+(cz)^2\otimes
    Y^+(q^2cz)\Psi(c^{3/2}z). \label{follows1}
  \end{eqnarray}
  
  We now calculate
  $$[\Delta(h_1),\sum_{n\geq 1}(-q(q-q^{-1})^2)^n
  X_0^+(cz)^{n+1}\otimes Y^+(q^2cz)^n\Psi(c^{3/2}z)].$$
  By Lemma \ref{lemma4} (\ref{lemma4.1}), we have that
  \begin{eqnarray*}
    \lefteqn{[h_1\otimes c^{3/2}, \sum_{n\geq 1}(-q(q-q^{-1})^2)^n
      X_0^+(cz)^{n+1}\otimes Y^+(q^2cz)^n\Psi(c^{3/2}z)]} \\
    &=& [2]z(c^{1/2}\otimes c^{3/2}) \\
    && \times \Big( \sum_{n\geq
      1}(-q(q-q^{-1})^2)^n(\sum_{k=0}^{n}q^{-2k}) X_0^+(cz)^{n+1}
    \otimes Y^+(q^2cz)^n \Psi(c^{3/2}z) \\
    && - \sum_{n\geq
      1}(-q(q-q^{-1})^2)^n(\sum_{k=0}^{n}q^{-2k}) x_0X_0^+(cz)^{n}
    \otimes Y^+(q^2cz)^n \Psi(c^{3/2}z)\Big),
  \end{eqnarray*}
  and, by Lemma \ref{lemma4} (\ref{lemma4.2}), that
  \begin{eqnarray*}
    \lefteqn{[c^{1/2}\otimes h_1, \sum_{n\geq 1}(-q(q-q^{-1})^2)^n
      X_0^+(cz)^{n+1}\otimes Y^+(q^2cz)^n\Psi(c^{3/2}z)]} \\
    &=& [2](c^{1/2}\otimes c^{1/2}) \\
    && \times \sum_{n\geq 1}(-q(q-q^{-1})^2)^n
    (\sum_{k=0}^{n-1}q^{-2k}) X_0^+(cz)^{n+1}\otimes
    Y^+(q^2cz)^{n-1}y_1 \Psi(c^{3/2}z) \\
    && -[2]z(c^{1/2}\otimes c^{3/2}) \\
    && \times \sum_{n\geq 1}(-q(q-q^{-1})^2)^n
    q^2(\sum_{k=0}^{n-1}q^{-2k})X_0^+(cz)^{n+1}\otimes
    Y^+(q^2cz)^{n}\Psi(c^{3/2}z).
  \end{eqnarray*}
  Further, by Lemma \ref{lemma4} (\ref{lemma4.5}), (\ref{lemma4.3})
  and (\ref{lemma4.8}),
  \begin{eqnarray*}
    \lefteqn{[x_0\otimes y_1, \sum_{n\geq 1}(-q(q-q^{-1})^2)^n
      X_0^+(cz)^{n+1}\otimes Y^+(q^2cz)^n\Psi(c^{3/2}z)]} \\
    &=& \sum_{n\geq 1}(-q(q-q^{-1})^2)^nx_0 X_0^+(cz)^{n+1} \\
    && \otimes \Big(
    q^{-2n}Y^+(q^2cz)^ny_1+(1-q^{-2n})q^2zcY^+(q^2cz)^{n+1}\Big)
    \Psi(c^{3/2}z) \\
    && -\sum_{n\geq 1}(-q(q-q^{-1})^2)^n\Big( q^{-2n-2}x_0
    X_0^+(cz)^{n+1} +(1-q^{-2n-2})X_0^+(cz)^{n+2}\Big) \\
    && \otimes Y^+(q^2cz)^{n}\Big( q^2y_1\Psi(c^{3/2}z)
    +(1-q^4)zcY^+(q^2cz)\Psi(c^{3/2}z)\Big).
  \end{eqnarray*}
  We will write
  \begin{eqnarray}
    \lefteqn{[\Delta(h_1),\sum_{n\geq 1}(-q(q-q^{-1})^2)^n
      X_0^+(cz)^{n+1}\otimes Y^+(q^2cz)^n\Psi(c^{3/2}z)]} \nonumber \\
    &=& [2]z(c^{1/2}\otimes c^{3/2})(F(z)+G(z))+[2](c^{1/2}\otimes
    c^{1/2})H(z), \label{follows2}
  \end{eqnarray}
  where $F(z),~G(z)$ and $H(z)$ are defined by
  \begin{eqnarray*}
    F(z) &=& \sum_{n\geq
      1}(-q(q-q^{-1})^2)^n(\sum_{k=0}^{n}q^{-2k}) X_0^+(cz)^{n+1}
    \otimes Y^+(q^2cz)^n \Psi(c^{3/2}z) \\
    && -\sum_{n\geq 1}(-q(q-q^{-1})^2)^n
    q^2(\sum_{k=0}^{n-1}q^{-2k})X_0^+(cz)^{n+1}\otimes
    Y^+(q^2cz)^{n}\Psi(c^{3/2}z) \\
    && -\sum_{n\geq 1}(-q(q-q^{-1})^2)^n
    (q-q^{-1})^2(q+q^{-1})q^2(1-q^{-2n-2}) \\
    && ~~X_0^+(cz)^{n+2}\otimes Y^+(q^2cz)^{n+1}\Psi(c^{3/2}z), \\
    G(z) &=& -\sum_{n\geq 1}(-q(q-q^{-1})^2)^n (\sum_{k=0}^{n}q^{-2k})
    x_0X_0^+(cz)^{n}\otimes Y^+(q^2cz)^{n}\Psi(c^{3/2}z) \\
    && -\sum_{n\geq 1}(-q(q-q^{-1})^2)^n(1-q^{-2n})q^2(q-q^{-1}) \\
    && ~~x_0X_0^+(cz)^{n+1}\otimes Y^+(q^2cz)^{n+1}\Psi(c^{3/2}z) \\
    && +\sum_{n\geq 1}(-q(q-q^{-1})^2)^n(q-q^{-1})(1-q^4)q^{-2n-2} \\
    && ~~x_0X_0^+(cz)^{n+1}\otimes Y^+(q^2cz)^{n+1}\Psi(c^{3/2}z), \\
    H(z) &=& \sum_{n\geq 1}(-q(q-q^{-1})^2)^n(\sum_{k=0}^{n-1}q^{-2k})
    X_0^+(cz)^{n+1}\otimes Y^+(q^2cz)^{n-1}y_1\Psi(c^{3/2}z) \\
    && +\sum_{n\geq 1}(-q(q-q^{-1})^2)^n(1-q^{-2n-2})q^2(q-q^{-1}) \\
    && ~~X_0^+(cz)^{n+2}\otimes Y^+(q^2cz)^{n}y_1\Psi(c^{3/2}z). 
  \end{eqnarray*}
  We consider the series $F(z), ~G(z)$ and $H(z)$ separately.

  First,
  \begin{eqnarray*}
    F(z) &=& -q(q-q^{-1})^2(1+q^{-2}-q^2)X_0^+(cz)^{2}\otimes
    Y^+(q^2cz)\Psi(c^{3/2}z) \\
    && +\sum_{n\geq 2}(-q(q-q^{-1})^2)^n \Big( \sum_{k=0}^{n}q^{-2k}
    -\sum_{k=-1}^{n-2}q^{-2k} +(q^2+1)(1-q^{-2n})\Big)  \\
    && ~~X_0^+(cz)^{n+1}\otimes Y^+(q^2cz)^{n}\Psi(c^{3/2}z) \\
    &=& (q-q^{-1})^3(1+q^2)X_0^+(cz)^{2}\otimes
    Y^+(q^2cz)\Psi(c^{3/2}z) \\
    &&  +\sum_{n\geq 1}(-q(q-q^{-1})^2)^nX_0^+(cz)^{n+1}\otimes
    Y^+(q^2cz)^{n}\Psi(c^{3/2}z). 
  \end{eqnarray*}
  We then have,
  \begin{eqnarray*}
    G(z) &=& -(1+q^{-2})(-q)(q-q^{-1})^2x_0X_0^+(cz)\otimes
    Y^+(q^2cz)\Psi(c^{3/2}z) \\
    && +\sum_{n\geq 2}(-q(q-q^{-1})^2)^n \Big( -\sum_{k=0}^{n}q^{-2k}
    +\sum_{k=0}^{n-2}q^{-2k} +(1+q^2)q^{-2n}\Big) \\
    && ~~x_0X_0^+(cz)^{n}\otimes Y^+(q^2cz)^{n}\Psi(c^{3/2}z) \\
    &=& (q-q^{-1})^2(q+q^{-1})x_0X_0^+(cz)\otimes
    Y^+(q^2cz)\Psi(c^{3/2}z).
  \end{eqnarray*}
  We finally have that
  \begin{eqnarray*}
    H(z) &=& -q(q-q^{-1})^2X_0^+(cz)^2\otimes y_1\Psi(c^{3/2}z) \\
    && +\sum_{n\geq 2}(-q(q-q^{-1})^2)^n \Big(
    \sum_{k=0}^{n-1}q^{-2k} -\sum_{k=0}^{n-1}q^{-2k}\Big) \\
    && ~~X_0^+(cz)^{n+1}\otimes Y^+(q^2cz)^{n-1}y_1\Psi(c^{3/2}z) \\
    &=& -q(q-q^{-1})^2X_0^+(cz)^2\otimes y_1\Psi(c^{3/2}z).
  \end{eqnarray*}
  Now combining (\ref{follows1}) and (\ref{follows2}) with the
  expressions $F(z), ~G(z)$ and
  $H(z)$ above, we get that
  \begin{eqnarray*}
    \lefteqn{[\Delta(h_1),1\otimes X_0^+(c^2z)+X_0^+(cz)\otimes
      \Psi(c^{3/2}z)} \\
    \lefteqn{~~ +\sum_{n\geq 1}(-q(q-q^{-1})^2)^n
      X_0^+(cz)^{n+1}\otimes Y^+(q^2cz)^n\Psi(c^{3/2}z)]} \\
    &=& [2](zc^{1/2}\otimes zc^{3/2}) \\
    && \Big( 1\otimes X_0^+(c^2z)+X_0^+(cz)\otimes
    \Psi(c^{3/2}z) \\
    && +\sum_{n\geq 1}(-q(q-q^{-1})^2)^n
    X_0^+(cz)^{n+1}\otimes Y^+(q^2cz)^n\Psi(c^{3/2}z) \\
    && -x_0\otimes K-1\otimes x_0\Big),
  \end{eqnarray*}
  which proves (\ref{thm5.1}).

  Now, the relation (\ref{thm5.4}) follows by using
  $$ \alpha\otimes \alpha \circ \tau \circ \Delta = \Delta \circ\alpha,$$
  {}from Proposition \ref{proposition3}. The relation (\ref{thm5.3}) is
  proved as (\ref{thm5.1}), by showing that both sides of
  (\ref{thm5.3}) is a solution of the functional equation
  $$[F(z),\Delta(h_1)]=[2](zc^{1/2}\otimes
  c^{3/2})F(z)-[2](c^{1/2}\otimes c^{1/2})\Delta(y_1),$$
  where
  $$F(z)=\sum_{k\geq 1}f_kz^{-k}, ~f_1=\Delta(y_1)(1\otimes c^{-1}),$$
  obtained from the relation
  $$[y_k,h_1]=[2]c^{1/2}y_{k+1}.$$
  Then the relation (\ref{thm5.2}) follows from (\ref{thm5.3}) by
  using
  $$ \alpha\otimes \alpha \circ \tau \circ \Delta = \Delta \circ\alpha.$$

  We now prove (\ref{thm5.5}). By Proposition \ref{proposition1} we
  have
  \begin{eqnarray*}
    \lefteqn{\sum_{k\geq 0}\Delta(\psi_k)(zc^{1/2}\otimes
      c^{3/2})^{-k}} \\
    &=& (q-q^{-1}) \sum_{k\geq 0}\Delta([x_k,y_0]c^{-k/2})(zc^{1/2}\otimes
    c^{3/2})^{-k} +\Delta(K^{-1}) \\
    &=& (q-q^{-1}) [\sum_{k\geq 0}\Delta(x_k)(zc\otimes
    c^{2})^{-k},\Delta(y_0)] +\Delta(K^{-1}).
  \end{eqnarray*}
  We calculate
  $$[\sum_{k\geq 0}\Delta(x_k)(zc\otimes c^{2})^{-k},y_0\otimes
  1+K^{-1}\otimes y_0],$$
  with
  \begin{eqnarray*}
    \lefteqn{\sum_{k\geq 0} \Delta(x_k)(zc\otimes c^2)^{-k}} \\
    &=& 1\otimes X_0^+(c^2z)+ \sum_{n\geq 0}(-q(q-q^{-1})^2)^n
    X_0^+(cz)^{n+1}\otimes Y^+(q^2cz)^n\Psi(c^{3/2}z).
  \end{eqnarray*}
  Now,
  $$[1\otimes X_0^+(c^2z),y_0\otimes 1]=0,$$
  and
  \begin{eqnarray*}
    \lefteqn{[\sum_{n\geq 0}(-q(q-q^{-1})^2)^n
    X_0^+(cz)^{n+1}\otimes Y^+(q^2cz)^n\Psi(c^{3/2}z),y_0\otimes 1]}
    \\
    &=& \sum_{n\geq 0}(-q(q-q^{-1})^2)^n (q-q^{-1})^{-1}q^n[n+1] \\
    &&(q^{-2n}\Psi(c^{1/2}z)X_0^+(q^2cz)^{n}-K^{-1}X_0^+(cz)^{n})\otimes
    Y^+(q^2cz)^{n}\Psi(c^{3/2}z),
  \end{eqnarray*}
  by Lemma \ref{lemma4} (\ref{lemma4.9}). We also have
  $$[1\otimes X_0^+(c^2z),K^{-1}\otimes
  y_0]=(q-q^{-1})^{-1}(K^{-1}\otimes \Psi(c^{3/2}z)-K^{-1}\otimes
  K^{-1}),$$
  and
  \begin{eqnarray*}
    \lefteqn{[X_0^+(cz)\otimes \Psi(c^{3/2}z),K^{-1}\otimes y_0]} \\
    &=& q^2K^{-1}X_0^+(cz)\otimes
    \Big( q^2y_0\Psi(c^{3/2}z)-(q^2-q^{-2})Y_0^+(q^2cz)\Psi(c^{3/2}z)\Big) \\
    && -K^{-1}X_0^+(cz)\otimes y_0\Psi(c^{3/2}z) \\
    &=& -q^2(q^2-q^{-2})K^{-1}X_0^+(cz)\otimes
    Y^+(q^2cz)\Psi(c^{3/2}z),
  \end{eqnarray*}
  where we used Lemma \ref{lemma4} (\ref{lemma4.7}) after applying
  $\beta$ on it. We also used the fact that
  $$Y_0^+(q^2cz)=y_0+Y^+(q^2cz).$$
  Finally, if we apply $\beta$ on Lemma \ref{lemma4} (\ref{lemma4.3}),
  and then use it and Lemma \ref{lemma4} (\ref{lemma4.6}), we get
  \begin{eqnarray*}
    \lefteqn{[\sum_{n\geq 1}(-q(q-q^{-1})^2)^nX_0^+(cz)^{n+1}\otimes
      Y^+(q^2cz)^{n}\Psi(c^{3/2}z), K^{-1}\otimes y_0]} \\
    &=& \sum_{n\geq 1}(-q(q-q^{-1})^2)^nq^{2n+2}K^{-1}X_0^+(cz)^{n+1}
    \otimes Y^+(q^2cz)^{n} \\
    && ~~\Big(
    q^2y_0\Psi(c^{3/2}z)-(q^2-q^{-2})Y_0^+(q^2cz)\Psi(c^{3/2}z)\Big)
    \\
    && -\sum_{n\geq 1}(-q(q-q^{-1})^2)^nK^{-1}X_0^+(cz)^{n+1} \\
    && \otimes \Big(
    q^{2n}Y^+(q^2cz)^ny_0-(1-q^{2n})Y^+(q^2cz)^{n+1}\Big)
    \Psi(c^{3/2}z).
  \end{eqnarray*}
  Summing up, we get that
  \begin{eqnarray*}
    \lefteqn{[1\otimes X_0^+(c^2z)+\sum_{n\geq 0}(-q(q-q^{-1})^2)^n
      X_0^+(cz)^{n+1}\otimes 
      Y^+(q^2cz)^{n}\Psi(c^{3/2}z),} \\
    \lefteqn{y_0\otimes 1+K^{-1}\otimes y_0]} \\
    &=& \sum_{n\geq 0}(-1)^n(q-q^{-1})^{2n}(q-q^{-1})^{-1}[n+1] \\
    && ~~\Psi(c^{1/2}z)X_0^+(q^2cz)^{n}\otimes
    Y^+(q^2cz)^{n}\Psi(c^{3/2}z) \\
    && -\sum_{n\geq 2}(-q(q-q^{-1})^2)^n(q-q^{-1})^{-1}q^n[n+1] \\
    && ~~K^{-1}X_0^+(cz)^{n}\otimes Y^+(q^2cz)^{n}\Psi(c^{3/2}z) \\
    && -\sum_{n\geq 1}(-q(q-q^{-1})^2)^n\Big(
    q^{2n+2}(q^2-q^{-2})-(1-q^{2n})\Big) \\
    && ~~K^{-1}X_0^+(cz)^{n+1}\otimes Y^+(q^2cz)^{n+1}\Psi(c^{3/2}z) \\
    && -(q-q^{-1})^{-1}K^{-1}\otimes K^{-1} \\
    &=& (q-q^{-1})^{-1}\sum_{n\geq 0}(-1)^n(q-q^{-1})^{2n}[n+1] \\
    && ~~\Psi(c^{1/2}z)X_0^+(q^2cz)^{n}\otimes
    Y^+(q^2cz)^{n}\Psi(c^{3/2}z) \\
    && -\sum_{n\geq 2}\Bigg( (-q(q-q^{-1})^2)^n(q-q^{-1})^{-1}q^n[n+1] \\
    && +(-q(q-q^{-1})^2)^{n-1}\Big( q^{2n}(q^2-q^{-2})-(1-q^{2n-2})\Big) \Bigg) \\
    && ~~K^{-1}X_0^+(cz)^{n}\otimes Y^+(q^2cz)^{n}\Psi(c^{3/2}z) \\
    && -(q-q^{-1})^{-1}K^{-1}\otimes K^{-1} \\
    &=& (q-q^{-1})^{-1}\sum_{n\geq 0}(-1)^n(q-q^{-1})^{2n}[n+1] \\
    && ~~\Psi(c^{1/2}z)X_0^+(q^2cz)^{n}\otimes
    Y^+(q^2cz)^{n}\Psi(c^{3/2}z) \\
    && -(q-q^{-1})^{-1}K^{-1}\otimes K^{-1},
  \end{eqnarray*}
  since
  \begin{eqnarray*}
    \lefteqn{(-q(q-q^{-1})^2)^n(q-q^{-1})^{-1}q^n[n+1]} \\
    \lefteqn{~~+(-q(q-q^{-1})^2)^{n-1}\Big(
      q^{2n}(q^2-q^{-2})-(1-q^{2n-2})\Big) }\\
    &=& (-q(q-q^{-1})^2)^{n-1}\Big( -q^{n+1}(q^{n+1}-q^{-n-1}) \\
    && +q^{2n+2}-q^{2n-2}-1+q^{2n-2}\Big) =0.
  \end{eqnarray*}
  Now, we have that
  $$\Delta(K^{-1})=K^{-1}\otimes K^{-1},$$
  so (\ref{thm5.5}) is proved. The relation (\ref{thm5.6}) follows by
  applying
  $$ \alpha\otimes \alpha \circ \tau \circ \Delta = \Delta \circ\alpha,$$
  on (\ref{thm5.5}).

  The theorem is proved.
\end{proof}

\section{Powers of $X_0^+(z)$ and $Y_0^+(z)$}

In this final section we first study powers of the generating
functions, and then we find the formulas for the comultiplication of
the generators in $\hat{U}_q.$

For any positive integer $n$ we want to express the generating functions
$X_0^+(z)$ and $Y_0^+(z)$ as
\begin{eqnarray}
  X_0^+(z)^n &=& \sum_{0\leq m_n\leq \cdots\leq m_1}
  f_{m_n,\ldots,m_1}(z)x_{m_n}\cdots x_{m_1}, \label{form1} \\
  Y_0^+(z)^n &=& \sum_{0\leq m_n\leq \cdots\leq m_1}
  g_{m_n,\ldots,m_1}(z)y_{m_1}\cdots y_{m_n},\label{form2} 
\end{eqnarray}  
where $f_{m_n,\ldots,m_1}(z),~g_{m_n,\ldots,m_1}(z)~\in \mathbb{C}(q)(z)$. Actually, we will
find that 
$$f_{m_n,\ldots,m_1}(z)=g_{m_n,\ldots,m_1}(z)\in \mathbb{Z}[q,q^{-1}][z^{-1}].$$

It is clear, by (\ref{prop1.8}) and (\ref{prop1.9}), that the functions
$f_{m_n,\ldots,m_1}$ and $g_{m_n,\ldots,m_1}$ exist, since if
$k>l+1>0$, (\ref{prop1.8}) can be written as
$$x_kx_l=q^2x_lx_k+q^2x_{k-1}x_{l+1}-x_{l+1}x_{k-1},$$
and if $k=l+1>0$, (\ref{prop1.8}) simplifies to
$$x_{l+1}x_l=q^2x_lx_{l+1},$$
and (\ref{prop1.9}) give similar rules for the $y_i$'s. \\
So commuting
$x_k$ and $x_l$ will always result in a finite sum of products
$x_ix_j$ with $0\leq i\leq j$.

It is also clear that
$$f_{m_n,\ldots,m_1}(z)=g_{m_n,\ldots,m_1}(z),$$
since the formula (\ref{form2}) follows from applying the
antiautomorphism $\beta$ from Proposition \ref{proposition3} on
(\ref{form1}).

Given an $n$-tuple $(m_n,\ldots,m_1),$ with $0\leq m_n\leq \cdots\leq m_1$,
we can define a $j$-tuple $(l_1,\ldots,l_j)$ as follows:\\
$l_1$ is the number of occurences of $m_n$ in $(m_n,\ldots,m_1)$;\\
$l_2$ is the number of occurences of $m_{n-l_1}$ in
$(m_n,\ldots,m_1)$;\\ \nopagebreak
\vdots\\ \nopagebreak[4]
$l_j$ is the number of occurences of $m_{n-l_1-l_2-\cdots -l_{j-1}}$ in
$(m_n,\ldots,m_1)$.\\
So $n=\sum_{i=1}^{j} l_i$.

{\em Example.} The $8$-tuple $(1,1,2,2,3,5,5,5)$ defines the $4$-tuple
$(2,2,1,3)$.\\

In the proof of Theorem \ref{theorem6} we will need the following mappings.\\
Define, for $(m_k,\ldots,m_1)$ as above, $c_{m_k,\ldots,m_1}(q)=1$ when $k=0,$
$$c_{m_1}(q)=1,$$
and inductively,
$$c_{m_k,\ldots,m_1}(q)=q^{k(k-1)m_k}\sum_{j=0}^{k-1}
(-1)^{k+j+1}\left[ {k \atop j} \right] q^{(k-j)(k-1)-j(j-1)m_k}
c_{m_j,\ldots,m_1}(q)\delta_{m_{j+1},m_k}.$$
Note that
\begin{equation}
  c_{m_k-a,\ldots,m_1-a}(q)=q^{-k(k-1)a}c_{m_k,\ldots,m_1}(q),
  \label{twostar}
\end{equation}
(by induction) for $0\leq a \leq m_k.$
\\

\begin{lemma} \label{lemma9}
  We have that
  \begin{equation}
    c_{m_n,\ldots,m_1}(q)=\frac{[n]!}{[l_1]!\cdots [l_j]!} 
    q^{\sum_{i=1}^{n} 2(i-1)m_i+\big( n(n-1)-\sum_{i=1}^{j}
      l_i(l_i-1)\big) /2}. \label{lemma9.1}
  \end{equation}
\end{lemma}

\begin{proof}
  For the proof we need the well-known identity (cf \cite{Guide},
  \cite{Lu})
  \begin{equation}
    \sum_{k=0}^{n} (-1)^{n+k} \left[ {n \atop k} \right]
    q^{(n-k)(n-1)} = 0, \label{onestar}
  \end{equation}
  for $n>0.$ \\
  We prove the lemma by induction over $n, ~n \geq 1.$ By definition,
  $$1=c_{m_1}(q)= q^0 \frac{[1]!}{[1]!}q^0.$$
  Suppose that the lemma is true for all $k$-tuples
  $(m_k,\ldots,m_1), ~0\leq m_k\leq\cdots\leq m_1,$ with $1 \leq k\leq n-1.$
  Fix $(m_n,\ldots,m_1)$ and the corresponding $(l_1,\ldots,l_j).$ We
  have, by definition,
  $$c_{m_n,\ldots,m_1}(q)=
  q^{n(n-1)m_n}\sum_{k=0}^{n-1} 
  (-1)^{n+k+1}\left[ {n \atop k} \right] q^{(n-k)(n-1)-k(k-1)m_n}
  c_{m_k,\ldots,m_1}(q)\delta_{m_{k+1},m_n}.$$
  Let $k$ be an integer such that $n-l_1\leq k\leq n-1.$ By the
  induction hypothesis,
  \begin{eqnarray*}
    c_{m_k,\ldots,m_1}(q) &=& \frac{[k]!}{[l_1-n+k]!\cdots [l_j]!} \\
    &&\times q^{\sum_{i=1}^{k}
      2(i-1)m_i+\big( k(k-1)-(l_1-n+k)(l_1-n+k-1)-\sum_{i=2}^{j}
      l_i(l_i-1)\big) /2},
  \end{eqnarray*}
  so we get that
  \begin{eqnarray*}
    c_{m_n,\ldots,m_1}(q) &=& \sum_{k=n-l_1}^{n-1} (-1)^{n+k+1}
    \left[{n \atop k} \right] \frac{[k]!}{[l_1-n+k]!\cdots [l_j]!} \\
    && \times q^{n(n-1)m_n+\sum_{i=1}^{k}2(i-1)m_i-k(k-1)m_n} \\
    && \times q^{(n-k)(n-1)+\big( k(k-1)-(l_1-n+k)(l_1-n+k-1)-
      \sum_{i=2}^{j}l_i(l_i-1)\big) /2}.
  \end{eqnarray*}
  We have that
  $m_{n-l_1+1}=\cdots=m_n,$ so since
  \begin{eqnarray*}
    \sum_{i=1}^{n} 2(i-1)m_i -\sum_{i=1}^{k} 2(i-1)m_i &=&
    2\sum_{i=k+1}^{n} (i-1)m_i \\
    &=& 2m_n\sum_{i=k+1}^{n} (i-1) \\
    &=& 2m_n(n-k)\frac{(n-1)+k}{2}  \\
    &=& m_n(n^2-k^2-n+k) \\
    &=& m_n n(n-1)-m_n k(k-1),
  \end{eqnarray*}
  we get that
  $$n(n-1)m_n+\sum_{i=1}^{k}2(i-1)m_i-k(k-1)m_n =
  \sum_{i=1}^{n}2(i-1)m_i.$$

  We now consider the expression
  $$(n-k)(n-1)+\Big( k(k-1)-(l_1-n+k)(l_1-n+k-1)\Big) /2.$$
  We have

  \lefteqn{(n-k)(n-1)+\frac{k(k-1)}{2}-\frac{l_1(l_1-1-(n-k))}{2}+\frac{(n-k)(l_1-1-(n-k))}{2}} 
  
  $$~=n(n-1)-k(n-1)+\frac{k(k-1)}{2}-\frac{l_1(l_1-1)}{2}+l_1(n-k)$$

  \begin{eqnarray*}
    -\frac{n-k}{2}-\frac{(n-k)^2}{2} &=& \frac{n(n-1)}{2}+(l_1-1)(n-k) 
    -\frac{l_1(l_1-1)}{2}.
  \end{eqnarray*}

  So, if we introduce the notation
  \begin{eqnarray*}
    A &=& q^{\sum_{i=1}^{n}2(i-1)m_i}, \\
    B &=& q^{\big( n(n-1)-\sum_{i=1}^{j}l_i(l_i-1)\big) /2},
  \end{eqnarray*}
  we have that
  \begin{eqnarray*}
    c_{m_n,\ldots,m_1}(q) &=& AB \frac{[n]!}{[l_1]!\cdots[l_j]!} 
    \sum_{k=n-l_1}^{n-1}(-1)^{n+k+1}
    \frac{[l_1]!}{[n-k]![l_1-n+k]!} q^{(l_1-1)(n-k)}\\
    &=& AB \frac{[n]!}{[l_1]!\cdots[l_j]!} 
    \sum_{r=0}^{l_1-1}(-1)^{l_1+r+1}
    \frac{[l_1]!}{[l_1-r]![r]!} q^{(l_1-1)(l_1-r)}.
  \end{eqnarray*}
  By (\ref{onestar}), we finally get
  $$c_{m_n,\ldots,m_1}(q) = AB \frac{[n]!}{[l_1]!\cdots[l_j]!}.$$

  The lemma is proved.
\end{proof}

We now state the main theorem of this section.
  
\begin{theorem} \label{theorem6}
  For any integer $n\geq 1,$ we have
  \begin{eqnarray}
    X_0^+(z)^n &=& \sum_{0\leq m_n\leq \cdots\leq m_1}
    c_{m_n,\ldots,m_1}(q)x_{m_n}\cdots x_{m_1}z^{-(m_1+\cdots+m_n)},
    \label{thm6.1} \\
    Y_0^+(z)^n &=& \sum_{0\leq m_n\leq \cdots\leq m_1}
    c_{m_n,\ldots,m_1}(q)y_{m_1}\cdots y_{m_n}z^{-(m_1+\cdots+m_n)},
    \label{thm6.2} 
  \end{eqnarray}
  where 
  \begin{eqnarray}
    c_{m_n,\ldots,m_1}(q)=\frac{[n]!}{[l_1]!\cdots [l_j]!}
    q^{\sum_{i=1}^{n} 2(i-1)m_i+\big( n(n-1)-\sum_{i=1}^{j}
      l_i(l_i-1)\big) /2},
    \label{thm6.3}
  \end{eqnarray}
  is in $\mathbb{Z}[q,q^{-1}]$.
\end{theorem}

We will need the following lemma before we can prove Theorem
\ref{theorem6}. The lemma tells us how the automorphism
$S:\mathcal{X}\to\mathcal{X}$ from Proposition \ref{proposition2}
acts on $X_0^+(z)^n.$

\begin{lemma} \label{lemma7}
  We have, for any $n\geq0,$
  \begin{eqnarray}
    S(X_0^+(z)^n) &=& z^n \sum_{k=0}^{n}(-1)^k \left[ {n \atop k}\right]
    q^{-(n-k)(n-1)}x_0^kX_0^+(z)^{n-k} \label{lemma7.1} \\
    &=& z^n \sum_{k=0}^{n}(-1)^{n-k} \left[ {n \atop k}\right]
    q^{-k(n-1)}x_0^{n-k}X_0^+(z)^{k}. \label{lemma7.2}
  \end{eqnarray}
\end{lemma}

\begin{proof}
  We only prove (\ref{lemma7.1}), since (\ref{lemma7.2}) follows from
  (\ref{lemma7.1}) by changing summation order.

  For $n=1$ we have that
  $$S(X_0^+(z))=z(X_0^+(z)-x_0),$$
  so (\ref{lemma7.1}) is true in this case. Now suppose, by induction,
  that the lemma is true for some $n\geq 1.$ Then
  \begin{eqnarray*}
    S(X_0^+(z)^{n+1}) &=& S(X_0^+(z)^{n})S(X_0^+(z)) \\
    &=& z^{n+1}\Bigg( \sum_{k=0}^{n}(-1)^k\left[ {n \atop k}\right]
    q^{-(n-k)(n-1)}x_0^kX_0^+(z)^{n-k+1} \\
    && -\sum_{k=0}^{n}(-1)^k\left[ {n \atop k}\right]
    q^{-(n-k)(n-1)}\Big(q^{-2(n-k)}x_0^{k+1}X_0^+(z)^{n-k} \\
    && +(1-q^{-2(n-k)})x_0^kX_0^+(z)^{n-k+1}\Big)\Bigg),
  \end{eqnarray*}
  where we have used Lemma \ref{lemma4} (\ref{lemma4.3}). It follows
  that
  \begin{eqnarray*}
    \lefteqn{S(X_0^+(z)^{n+1})} \\
    &=& z^{n+1}\Big( \sum_{k=0}^{n}(-1)^{k+1}\left[
    {n \atop k}\right]
    q^{-(n-k)(n-1)}q^{-2(n-k)}x_0^{k+1}X_0^+(z)^{n-k} \\
    && +\sum_{k=0}^{n}(-1)^{k}\left[
    {n \atop k}\right]
    q^{-(n-k)(n-1)}q^{-2(n-k)}x_0^{k}X_0^+(z)^{n-k+1}\Big) \\
    &=& z^{n+1}\Bigg( \sum_{k=0}^{n+1}(-1)^{k}\Big( q^{k-(n+1)}\left[
    {n \atop k-1}\right] +q^k\left[ {n \atop k}\right] \Big)
    q^{-(n+1-k)n}x_0^kX_0^+(z)^{n+1-k}\Bigg) \\
    &=& z^{n+1}\sum_{k=0}^{n+1}(-1)^{k}\left[ {n+1 \atop k}\right]
    q^{-(n+1-k)n} x_0^kX_0^+(z)^{n+1-k}.
  \end{eqnarray*}
  The lemma is proved.
\end{proof}

We can now prove Theorem \ref{theorem6}.

\begin{proof}[Proof (of Theorem \ref{theorem6}).]
  It is enough to prove (\ref{thm6.1}), since (\ref{thm6.2}) follows
  {}from (\ref{thm6.1}) by applying the antiautomorphism $\beta$ from
  Proposition \ref{proposition3}.

  Let $n$ be a positive integer. Write
  $$X_0^+(z)^n=\sum_{0\leq m_n\leq \cdots\leq m_1} f_{m_n,\ldots,
    m_1}(z) x_{m_n}\cdots x_{m_1}.$$
  for some $f_{m_n,\ldots,m_1}(z)\in \mathbb{C}(q)(z).$ We will prove
  that
  \begin{equation}
    f_{m_n,\ldots,m_1}(z)
    =c_{m_n,\ldots,m_1}(q)z^{-(m_1+\cdots +m_n)}. \label{fzres}
  \end{equation}
  for $c_{m_n,\ldots,m_1}(q)$ as in (\ref{lemma9.1}). For $n=1$ it is
  clear that (\ref{fzres}) is true. Now suppose that $n>1,$ and that
  (\ref{fzres}) is true for all $k$ such that $1\leq k<n.$ By Lemma
  \ref{lemma7}, we have that
  \begin{eqnarray}
    \lefteqn{S(X_0^+(z)^n) = \sum_{1\leq m_n\leq \cdots\leq m_1}
    f_{m_n-1,\ldots,m_1-1}(z) x_{m_n}\cdots x_{m_1}}  \nonumber \\
    &=& z^n\sum_{k=0}^{n}(-1)^k\left[ {n \atop k}\right]
    q^{-(n-k)(n-1)} \sum_{0\leq s_{n-k}\leq \cdots\leq
      s_1}f_{s_{n-k},\ldots ,s_1}(z)x_0^k x_{s_{n-k}} \cdots
    x_{s_1}. \label{Sformula} 
  \end{eqnarray}

  Now fix $(m_n,\ldots,m_1)$ such that $1\leq m_n\leq \cdots\leq m_1.$
  By PBW for Quantum Affine algebras (see \cite{Be}), and by
  (\ref{Sformula}), we have that
  $$f_{m_n-1,\ldots,m_1-1}(z)=z^nq^{-n(n-1)}f_{m_n,\ldots,m_1}(z).$$
  It follows that
  \begin{equation}
    f_{m_n,\ldots,m_1}(z)=z^{-m_nn}q^{n(n-1)m_n}
    f_{0,m_{n-1}-m_n,\ldots,m_1-m_n}(z). \label{f0res}
  \end{equation}

  We now consider the terms where $x_{m_n}=x_0$ in
  (\ref{Sformula}). Extracting coefficients gives
  
  \begin{eqnarray*}
    0 &=& {\hbox{coeff. of\ }} x_0x_{m_{n-1}}\cdots x_{m_1}\hbox{ in} \Bigg(
    \left[ {n \atop 0}\right] q^{-n(n-1)} \sum_{0\leq s_{n}\leq \cdots\leq
      s_1}f_{s_{n},\ldots ,s_1}(z)x_{s_{n}} \cdots x_{s_1} \\
    && -\left[ {n \atop 1}\right] q^{-(n-1)(n-1)} \sum_{0\leq s_{n-1}\leq \cdots\leq
      s_1}f_{s_{n-1},\ldots ,s_1}(z)x_0x_{s_{n-1}} \cdots x_{s_1} \\
    && +\sum_{k=2}^{n}(-1)^k\left[ {n \atop k}\right]
    q^{-(n-k)(n-1)} \sum_{0\leq s_{n-k}\leq \cdots\leq
      s_1}f_{s_{n-k},\ldots ,s_1}(z)x_0^k x_{s_{n-k}} \cdots
    x_{s_1}\Bigg).
  \end{eqnarray*}
  When $m_{n-1}\geq 1,$ this gives
  $$q^{-n(n-1)}f_{0,m_{n-1},\ldots,m_1}(z) =\left[ {n \atop
    1}\right] q^{-(n-1)(n-1)}
  f_{m_{n-1},\ldots,m_1}(z),$$
  and when $m_{n-1}=0, ~m_{n-2}\geq 1,$ we get
  $$q^{-n(n-1)}f_{0,0,m_{n-2},\ldots,m_1}(z) $$
  $$ \quad =\left[ {n \atop 1}\right] q^{-(n-1)(n-1)}
  f_{0,m_{n-2},\ldots,m_1}(z)-\left[ {n \atop 2}\right]
  q^{-(n-2)(n-1)}f_{m_{n-2},\ldots,m_1}(z).$$
  It is easy to see that in general we have,
  $$q^{-n(n-1)}f_{0,m_{n-1},\ldots,m_1}(z) =
  \sum_{k=1}^{n}(-1)^{k+1}\left[ {n \atop k}\right] q^{-(n-k)(n-1)}
  f_{m_{n-k},\ldots,m_1}(z) \delta_{m_{n-k+1},0}.$$
  So, (\ref{f0res}) gives, by induction, and by (\ref{twostar}), that
  \begin{eqnarray*}
    \lefteqn{f_{m_n,\ldots,m_1}(z)} \\
    &=& q^{n(n-1)m_n} \sum_{k=0}^{n-1}(-1)^{n+1+k}\left[ {n
      \atop k}\right]
    q^{(n-k)(n-1)}c_{m_k-m_n,\ldots,m_1-m_n}(q)
    \delta_{m_{k+1},m_n}z^{-(m_1+\cdots +m_n)} \\
    &=& q^{n(n-1)m_n} \sum_{k=0}^{n-1}(-1)^{n+1+k}\left[ {n
      \atop k}\right]
    q^{(n-k)(n-1)-k(k-1)m_n}c_{m_k,\ldots,m_1}(q)
    \delta_{m_{k+1},m_n} z^{-(m_1+\cdots +m_n)} \\
    &=& c_{m_n,\ldots,m_1}(q)z^{-(m_1+\cdots +m_n)}.
  \end{eqnarray*}
  By Lemma \ref{lemma9}, we have proved the theorem.
\end{proof}  

We now want to find the explicit action of the comultiplication on the
generators of $\hat{U}_q.$ For this we define, for $n\geq 1, ~m\geq
0,$ the sets
$$\Omega_{n,m}=\{ (m_n,\ldots,m_1)\in \mathbb{Z}^n ; 0\leq m_n\leq
\cdots \leq m_1, ~\sum_{i=1}^{n}m_i=m\},$$
and
$$\Omega^+_{n,m}=\{ (m_n,\ldots,m_1)\in \mathbb{Z}^n ; 1\leq m_n\leq
\cdots \leq m_1, ~\sum_{i=1}^{n}m_i=m\}.$$
Then we have, by Theorem \ref{theorem6} and (\ref{twostar}),
\begin{eqnarray*}
  X_0^+(z)^n &=& \sum_{m\geq 0}\sum_{(a_n,\ldots, a_1)\in
    \Omega_{n,m}} c_{a_n,\ldots, a_1}(q)x_{a_n}\cdots x_{a_1} z^{-m},
  \\
  Y^+(z)^n &=& z^{-n}T(Y_0^+(z)^n) \\
  &=& \sum_{m\geq n}\sum_{(b_n,\ldots, b_1)\in
    \Omega^+_{n,m}} c_{b_n,\ldots, b_1}(q)q^{-n(n-1)}y_{b_1}\cdots y_{b_n} z^{-m},
\end{eqnarray*}
so that, by Theorem \ref{theorem5},
\begin{eqnarray}
  \lefteqn{\sum_{k\geq 0}\Delta (x_k)(zc\otimes c^2)^{-k}}  \nonumber \\
  &=& 1\otimes \sum_{k\geq 0}x_kc^{-2k}z^{-k}+ \sum_{k \geq
    0}x_kc^{-k}z^{-k} \otimes \sum_{l \geq 0}\psi_lc^{-3l/2}z^{-l}
  \nonumber \\
  &&+ \sum_{n\geq
    1}(-q(q-q^{-1})^2)^n \Big( \sum_{m\geq 0}\sum_{(a_{n+1},\ldots, a_1)\in
  \Omega_{n+1,m}} c_{a_{n+1},\ldots, a_1}(q)x_{a_{n+1}}\cdots
  x_{a_1} c^{-m}z^{-m}\Big) \nonumber \\
  && \otimes \Big( \sum_{m\geq n}\sum_{(b_n,\ldots, b_1)\in
    \Omega^+_{n,m}} c_{b_n,\ldots, b_1}(q)q^{-n(n-1)}y_{b_1}\cdots y_{b_n} q^{-2m}c^{-m}z^{-m}
  \Big) \nonumber \\
  && \times \Big( \sum_{m\geq 0}\psi_mc^{-3m/2}z^{-m} \Big).
  \label{storsum}
\end{eqnarray}
By identifying coefficients of $z^{-N}, ~N\geq 0,$ in (\ref{storsum}), 
and using the convention that summation over an empty set is zero,
we get that
\begin{eqnarray*}
  \Delta(x_N) &=& c^N\otimes x_N +\sum_{k=0}^{N}c^{N-k}x_k\otimes
    c^{(N+3k)/2}\psi_{N-k} \\
  && +\sum_{n=1}^{N}(-q(q-q^{-1})^2)^n q^{-n(n-1)} \\
  && \times \sum_{m=0}^{N}\Big( \sum_{(a_{n+1},\ldots, a_1)\in
  \Omega_{n+1,m}} c_{a_{n+1},\ldots, a_1}(q)c^{N-m}x_{a_{n+1}}\cdots
  x_{a_1} \\
  && \otimes \sum_{k=0}^{N-m}\sum_{(b_n,\ldots, b_1)\in
    \Omega^+_{n,N-m-k}} c_{b_n,\ldots,
    b_1}(q)q^{-2(N-m-k)}c^{N+m-k/2}y_{b_1}\cdots y_{b_n}\psi_k \Big).
\end{eqnarray*}
With similar reasoning we can find the comultiplication of the other
generators. We summarize the result in the following corollary.

\begin{corollary} \label{corollary7}
  For every integer $N\geq 1,$ the comultiplication satisfies \\
  \lefteqn{\Delta(x_N) = c^N\otimes x_N +\sum_{k=0}^{N}c^{N-k}x_k\otimes
    c^{(N+3k)/2}\psi_{N-k} }\\
  \lefteqn{~~+\sum_{n=1}^{N}(-q(q-q^{-1})^2)^n q^{-n(n-1)}} \\
    $$ \times \sum_{m=0}^{N}\Big( \sum_{(a_{n+1},\ldots, a_1)\in
    \Omega_{n+1,m}} c_{a_{n+1},\ldots, a_1}(q)c^{N-m}x_{a_{n+1}}\cdots
    x_{a_1} $$
    $$ \otimes \sum_{k=0}^{N-m}\sum_{(b_n,\ldots, b_1)\in
    \Omega^+_{n,N-m-k}} c_{b_n,\ldots,
    b_1}(q)q^{-2(N-m-k)}c^{N+m-k/2}y_{b_1}\cdots y_{b_n}\psi_k \Big); $$
  \lefteqn{\Delta(x_{-N}) = c^{-N}\otimes x_{-N} +\sum_{k=1}^{N}c^{-(N-k)}x_{-k}\otimes
    c^{(N-k)/2}\phi_{-(N-k)}} \\
  \lefteqn{~~+\sum_{n=1}^{N}(-q(q-q^{-1})^2)^n q^{n(n+1)}} \\ 
    $$ \times \sum_{m=0}^{N}\Big( \sum_{(a_{n+1},\ldots, a_1)\in
    \Omega^+_{n+1,m}} c_{a_{n+1},\ldots, a_1}(q^{-1})c^{-(N-m)}x_{-a_{n+1}}\cdots
    x_{-a_1}$$
    $$ \otimes \sum_{k=0}^{N-m}\sum_{(b_n,\ldots, b_1)\in
    \Omega_{n,N-m-k}} c_{b_n,\ldots,
    b_1}(q^{-1})q^{2(N-m-k)}c^{N-m-k/2}y_{-b_1}\cdots y_{-b_n}\phi_{-k}
    \Big); $$
  \lefteqn{\Delta(y_N) = y_N\otimes c^N +\sum_{k=1}^{N}c^{-(N-k)/2}\psi_{N-k}\otimes
    c^{N-k}y_k} \\
  \lefteqn{~~+\sum_{n=1}^{N}(-q^{-1}(q-q^{-1})^2)^nq^{-n(n+1)}} \\ 
    $$ \times \sum_{m=0}^{N}\Big( \sum_{k=0}^{N-m}\sum_{(a_{n},\ldots, a_1)\in
    \Omega_{n,N-m-k}} c_{a_{n},\ldots, a_1}(q)q^{-2(N-m-k)}c^{-(N-m-k/2)}\phi_{k}x_{a_{n}}\cdots
    x_{a_1} $$
    $$\qquad \otimes \sum_{(b_{n+1},\ldots, b_1)\in
    \Omega^+_{n+1,m}} c_{b_n,\ldots,
    b_1}(q)c^{N-m}y_{b_1}\cdots y_{b_{n+1}} \Big); $$
  \lefteqn{\Delta(y_{-N}) = y_{-N}\otimes c^{-N}
    +\sum_{k=0}^{N}c^{-(N+3k)/2}\phi_{-(N-k)} \otimes
    c^{-(N-k)}y_{-k}} \\
  \lefteqn{~~+\sum_{n=1}^{N}(-q^{-1}(q-q^{-1})^2)^n q^{n(n-1)}} \\
    $$ \times \sum_{m=0}^{N}\Big( \sum_{k=0}^{N-m}\sum_{(a_n,\ldots, a_1)\in
    \Omega^+_{n,N-m-k}} c_{a_n,\ldots,
    a_1}(q^{-1})q^{2(N-m-k)}c^{-(N+m-k/2)}\phi_{-k}x_{-a_n}\cdots
    x_{-a_1} $$
    $$ \otimes \sum_{(b_{n+1},\ldots, b_1)\in
    \Omega_{n+1,m}} c_{b_{n+1},\ldots, b_1}(q^{-1})c^{-(N-m)}y_{-b_{1}}\cdots
    y_{-b_{n+1}} \Big); $$
  \lefteqn{\Delta(\psi_N) = \sum_{k=0}^{N}c^{(N-k)/2}\psi_k\otimes
    c^{3k/2}\psi_{N-k}} \\
  \lefteqn{~~+\sum_{n=1}^{N}(-1)^n(q-q^{-1})^{2n}[n+1]q^{-n(n-1)}} \\
    $$ ~\times \sum_{m=0}^{N}\Big( \sum_{k=0}^{m}\sum_{(a_{n},\ldots, a_1)\in
    \Omega_{n,m-k}} c_{a_{n},\ldots,
    a_1}(q)q^{-2(m-k)}c^{N/2-m+k/2}\psi_kx_{a_{n}}\cdots
    x_{a_1} $$
    $$ ~\otimes \sum_{l=0}^{N-m}\sum_{(b_n,\ldots, b_1)\in
    \Omega^+_{n,N-m-l}} c_{b_n,\ldots,
    b_1}(q)q^{-2(N-m-l)}c^{N/2+m-l/2}y_{b_1}\cdots y_{b_n}\psi_l \Big); $$
  \lefteqn{\Delta(\phi_{-N}) =\sum_{k=0}^{N}c^{-3k/2}\phi_{-(N-k)}\otimes
    c^{-(N-k)/2}\phi_{-k}} \\
  \lefteqn{~~+\sum_{n=1}^{N}(-1)^n(q-q^{-1})^{2n}[n+1]q^{n(n-1)}} \\
    $$ \times \sum_{m=0}^{N}\Big( \sum_{l=0}^{N-m}\sum_{(a_{n},\ldots, a_1)\in
    \Omega^+_{n,N-m-l}} c_{a_{n},\ldots,
    a_1}(q^{-1})q^{2(N-m-l)}c^{-(N/2+m-l/2)}\phi_{-l}x_{-a_{n}}\cdots
    x_{-a_1} $$
    $$ \otimes \sum_{k=0}^{m}\sum_{(b_n,\ldots, b_1)\in
    \Omega_{n,m-k}} c_{b_n,\ldots,
    b_1}(q^{-1})q^{2(m-k)}c^{-(N/2-m+k/2)}y_{-b_1}\cdots
    y_{-b_n}\phi_{-k} \Big),$$
  where empty sums are zero.
\end{corollary}

Mathematics (Faculty of Science), \\
Centre for Mathematical Sciences, \\
Lund University, \\
Box 118, \\
221 00 LUND, Sweden. \\
{\em email:} {\bf thoren@maths.lth.se}

\end{document}